\documentclass[preprint,12pt]{amsart}

\usepackage[margin=1.5in]{geometry}

\usepackage{amsmath}
\usepackage{amsfonts}
\usepackage{amssymb}
\usepackage{xypic}
\usepackage{longtable}
\usepackage{amsthm}
\usepackage{enumerate}

\theoremstyle{plain}
\newtheorem{theorem}{Theorem}[section]
\newtheorem{corollary}[theorem]{Corollary}

\theoremstyle{definition}
\newtheorem{definition}[theorem]{Definition}

\newtheorem{rem}[theorem]{Remark}

\makeatletter
\def\ps@pprintTitle{%
  \let\@oddhead\@empty
  \let\@evenhead\@empty
  \let\@oddfoot\@empty
  \let\@evenfoot\@oddfoot
}
\makeatother

\title{A note on Galois embeddings of abelian varieties}
\author{Robert Auffarth}

\address{Departamento de Matem\'aticas, Facultad de
Ciencias, Universidad de Chile, Santiago\\Chile}
\email{rfauffar@uchile.cl}

\thanks{Partially supported by Fondecyt Grant 3150171 and CONICYT PIA ACT1415}

\subjclass[2010]{14K05; 14E20}
\keywords{abelian variety, Galois embedding, automorphism}

\begin{document}

\maketitle

\begin{abstract}
In this note we show that if an abelian variety possesses a Galois embedding into some projective space, then it must be isogenous to the self product of an elliptic curve. We prove moreover that the self product of an elliptic curve always has infinitely many Galois embeddings. \end{abstract}

\section{Introduction}

Let $X$ be a smooth irreducible projective variety of dimension $d$ over an algebraically closed field $k$ of characteristic zero and let $D$ be a very ample divisor on $X$ that induces an embedding $\varphi_D:X\hookrightarrow\mathbb{P}^N$. Let $W\in\mathbb{G}(N-d-1,N)$ be a linear subvariety of $\mathbb{P}^N$ that is disjoint from the image of $X$, and let $W_0\in\mathbb{G}(d,N)$ be a linear subvariety disjoint from $W$. Consider the linear projection
$$\pi_W:\mathbb{P}^N\dashrightarrow W_0$$
with center $W$ and its restriction $f_W=\pi_W\varphi_D$ to $X$ (which is regular). We will write $k(X)$ for the function field of $X$.

\begin{definition}
The embedding $\varphi_D$ is called \emph{Galois} if there exists a subspace $W\in\mathbb{G}(N-d-1,N)$ such that $f_W$ is a Galois cover; that is, the field extension $f_W^*k(\mathbb{P}^d)\subseteq k(X)$ is Galois. In this case $W$ will be called a \emph{Galois subspace} for $D$.
\end{definition}

A question studied in \cite{Yoshi} is whether a Galois embedding exists for a smooth projective variety. Yoshihara proves that if $X$ is smooth and projective and $D$ is very ample, then the pair $(X,D)$ defines a Galois embedding if and only if the following conditions hold:

\begin{itemize}
\item There exists a finite subgroup $G$ of $\mbox{Aut}(X)$ that preserves $D$ and such that $|G|=(D^d)$.
\item There exists a $G$-invariant linear subspace $V$ of $H^0(X,\mathcal{O}_X(D))$ of dimension $d+1$ such that for every $\sigma\in G$, $\sigma^*|_{\mathbb{P}V}=\mbox{id}$.
\item $V$ is base-point free. 
\end{itemize}

Essentially, this shows that one wants to find a finite subgroup $G\leq\mbox{Aut}(X)$ such that $X/G\simeq\mathbb{P}^d$ and $h^*\mathcal{O}_{\mathbb{P}^d}(1)$ is very ample, where $h:X\to\mathbb{P}^d$ is the natural projection.

If a very ample divisor does not have a Galois embedding, then one can ask about the Galois closure of the extension field described above. It is also interesting to study what kinds of groups appear as Galois groups of Galois closures associated to very ample divisors. These types of questions have been widely studied for plane curves by Miura \cite{Miura}, Pirola and Schlesinger \cite{PS}, Takahashi \cite{Taka} and Yoshihara \cite{Yoshi1}, among many others. Progress has also been made for surfaces and hypersurfaces (\cite{Yoshi2}, \cite{Yoshi3}, among others), but much less is known in general for high dimensional cases (see Cukierman \cite{Cuk} for a general study of these questions).

In \cite{Yoshi}, Yoshihara restricts the above question to when $X$ is an abelian variety, and shows in \cite[Corollary 3.2]{Yoshi} that if $X$ is an abelian variety that possesses a Galois embedding, it must contain an elliptic curve. This happens because each irreducible component of the ramification divisor of the projection is a translate of an abelian subvariety of codimension 1. It is then shown that if an abelian surface possesses a Galois embedding, it is isogenous to the self-product of an elliptic curve. 

Let $\sim$ denote the isogeny relation. Our result generalizes the behavior observed in abelian surfaces to arbitrary dimension:

\begin{theorem}\label{theo1}
Let $A$ be an abelian variety of dimension $d$. 

\begin{enumerate}
\item\label{1} If $A$ possesses a Galois embedding into some projective space, then there exists an elliptic curve $E$ such that $A\sim E^d$.
\item\label{2} If $A=E^d$, then $A$ possesses infinitely many Galois embeddings with Galois groups of arbitrarily large order.
\end{enumerate}
\end{theorem}

We will give two proofs of part $(2)$ of the above theorem that produce different Galois groups. Moreover, the first proof of part $(2)$ gives us the following easy corollary:

\begin{corollary}
Let $C$ be a smooth projective curve that possesses a Galois embedding with Galois group $G$. Then $C^d$ possesses a Galois embedding with Galois group $G^d\rtimes S_d$, where $S_d$ is the symmetric group on $d$ letters.
\end{corollary}
 
We remark a subtlety involved in our description of abelian varieties that possess Galois embeddings. It is not clear if an abelian variety that is only isogenous, and not isomorphic, to the self product of an elliptic curve always has a Galois embedding. Indeed, the automorphism group of an abelian variety can change under an isogeny. It would be interesting to analyze this situation. Note that part (\ref{2}) of the theorem produces automorphism groups of arbitrarily large order; this is only possible due to the fact that the Galois groups will contain translations. Indeed, the order of a group of automorphisms that fix the origin and preserve a given polarization on an abelian variety is bounded above by a constant that depends only on the dimension of the abelian variety.

Another interesting question is to ask about what the minimal order of a Galois group can be. The Galois groups found in this paper for $E^d$ are of orders $2^dd!|Q_0|^d$ and $(d+1)!|Q_0|^d$ where $Q_0$ is any non-trivial finite subgroup of $E$. Note that by the Riemann-Roch Theorem for abelian varieties, if a very ample divisor $D$ on $A$ induces a Galois embedding with Galois group $G$, then $|G|=(D^d)=d!\chi(\mathcal{O}_A(D))$, and so $d!$ divides $|G|$. We note that for abelian surfaces, however, the groups found in this paper do not give a minimal Galois embedding. Indeed, for $d=2$ the smallest group we have found is of order 24, but the main result of \cite{Yoshi4} implies that the minimal degree that appears for abelian surfaces is 16.\\

\noindent\textit{Acknowledgements:} I would like to thank Dr. Hisao Yoshihara for several beneficial conversations, as well as for generously sharing his research with me.

\section{Proof of Theorem \ref{theo1}}

\subsection{Proof of part (\ref{1})} Let $A$ be an abelian variety, and assume that it has a Galois embedding given by a group $G$. We have the exact sequence
$$1\to G_0\to G\to H\to 1$$
where $G_0$ are the translations that lie in $G$ and $H=\{f-f(0):f\in G\}$. Let $B$ be the abelian variety $A/G_0$; we notice that $A/G\simeq B/H$. Now by Yoshihara's results in \cite{Yoshi}, $A$ contains an elliptic curve. He proves this by showing that the reduced components of the ramification divisor are abelian subvarieties of codimension 1. Another way of showing the existence of an elliptic curve in $A$ is to note that since the quotient $A/G$ is smooth, $G$ must be generated by pseudoreflections. In particular, there exists a $\sigma\in G$ such that it fixes a divisor pointwise. However, the irreducible components of such a divisor correspond to certain irreducible components of $\ker(\mbox{id}-\sigma)$ which are (translations of) abelian subvarieties of codimension 1. By any means, $A$ contains an elliptic curve, and therefore $B$ also contains some elliptic curve that we will name $E$. If $\pi:B\to\mathbb{P}^d$ is the projection, set $\mathcal{L}:=\pi^*\mathcal{O}_{\mathbb{P}^d}(1)$.

Assume that $B\sim E^r\times Z$ where $Z$ is a $(d-r)$-dimensional abelian subvariety of $B$. Let $Y$ be a complementary abelian subvariety of $Z$ in $B$ with respect to $\mathcal{L}$, so that $Y\sim E^r$, and let $N_Y,N_Z\in\mbox{End}(B)$ be the norm endomorphisms associated to $Y$ and $Z$. Recall that
$$N_Z=e_Zi_Z\phi_{\mathcal{L|_Z}}^{-1}i_Z^\vee\phi_\mathcal{L},$$
where $i_Z:Z\hookrightarrow B$ is the inclusion, $\phi_\mathcal{L}:B\to B^\vee:=\mbox{Pic}^0(B)$ (respectively $\phi_{\mathcal{L}|_Z}:Z\to Z^\vee$) is the isogeny induced by $\mathcal{L}$ (respectively $\mathcal{L}|_Z$) and $e_Z$ is the exponent of the finite group $\ker\phi_{\mathcal{L}|_Z}$ (see \cite[Section 5.3]{BL}). Essentially $N_Z$ is a sort of orthogonal projection from $B$ to $Z$, and satisfies $N_Z^2=e_ZN_Z$.

Let us assume that $\mbox{Hom}(Y,Z)=0$. In other words, assume $Z$ contains no elliptic curve that is isogenous to $E$. This implies that every endomorphism of $Y\times Z$ can be written as $(x,y)\mapsto (a(x),b(y))$ where $a\in\mbox{End}(Y)$ and $b\in\mbox{End}(Z)$. In particular, every endomorphism preserves $\{0\}\times Z$. If $\sigma\in H$, then we have the commutative diagram\\

\centerline{\xymatrix{B\ar[d]_\sigma&&Y\times Z\ar[ll]_-{N_Y+N_Z}\ar[d]^{\tilde{\sigma}}\\
B\ar[rr]^-{(N_Y,N_Z)}&&Y\times Z}}

\vspace{0.3cm}

\noindent where $\tilde{\sigma}$ is simply the composition of the other three morphisms. We see that by the previous discussion, $\sigma(Z)=Z$ and $\sigma(Y)=Y$.  

We notice that for every $\sigma\in H$, $N_Z$ commutes with $\sigma$. Indeed, if $x=y+z$ with $y\in Y$ and $z\in Z$, then
$$N_Z(\sigma(y+z))=N_Z(\sigma(z))=e_Z\sigma(z)=\sigma(e_Zz)=\sigma(N_Z(y+z)).$$
Therefore $N_Z$ descends to a morphism $\Phi_Z:\mathbb{P}^d\to\mathbb{P}^d$. Now
$$\mathcal{M}:=N_Z^*\mathcal{L}\simeq\pi^*\Phi_Z^*\mathcal{O}_{\mathbb{P}^d}(1)\simeq\pi^*\mathcal{O}_{\mathbb{P}^d}(k)$$
for some $k\in\mathbb{Z}$. By \cite[Proposition 3.1]{ALR}, 
$$(\mathcal{M}^d)=0\hspace{0.5cm}\mbox{and}\hspace{0.5cm}(\mathcal{M}\cdot\mathcal{L}^{d-1})\neq0.$$
But $(\mathcal{M}^d)=(\deg\pi)(\mathcal{O}_{\mathbb{P}^d}(k)^d)=|H|k^d$ which is 0 if and only if $k=0$. But if $k=0$, we have that $\mathcal{M}\simeq\mathcal{O}_B$, a contradiction.

Therefore, we must have that $\mbox{Hom}(Y,Z)\neq0$, and so $Z$ contains an elliptic curve that is isogenous to $E$. Therefore, $B\sim E^{r+1}\times Z'$ for some $Z'$, and we can continue the process until we have shown that $B\sim E^d$. \qed

\begin{rem}
We note that essentially the same proof shows that if $B=A/G_0$ is as above, then the action of $H$ on $B$ does not restrict to an action on a non-trivial abelian subvariety. Therefore, if $E\subseteq B$ is an elliptic curve, then $\sum_{\sigma\in H}\sigma(E)=B$. This shows that $B$ actually contains $d$ elliptic curves isomorphic (and not only isogenous) to $E$ whose pairwise intersection is finite. It is not clear, however, if $B$ is isomorphic to the product $E^d$.
\end{rem}

\subsection{First proof of part (\ref{2})} Let $C$ be a smooth projective curve and assume it has a Galois embedding induced by a group $G\leq\mbox{Aut}(C)$. Notice that $G^d\rtimes S_{d}$ acts on $C^d$: $G^d$ acts factor to factor and $S_d$ permutes the factors. It is also clear that the quotient of $C^d$ by $G^d\rtimes S_d$ gives $\mbox{Sym}^d(\mathbb{P}^1)\simeq\mathbb{P}^d$. Moreover, if $\mathcal{L}$ is the very ample line bundle that gives the embedding on $C$, then the pullback of $\mathcal{O}_{\mathbb{P}^d}(1)$ to $C^d$ is just $\mathcal{L}^{\boxtimes d}$ which is very ample. Therefore, if $C$ has a Galois embedding, so does $C^d$.

Now let $E$ be an elliptic curve. If $Q_0$ is a finite subgroup of $E$ (which we can see as a subgroup of translations of $\mbox{Aut}(E)$), we have the group of automorphisms $G=Q_0\rtimes\mathbb{Z}/2\mathbb{Z}$ where $\mathbb{Z}/2\mathbb{Z}$ acts on $E$ by $z\mapsto -z$. The quotient of $E$ by $G$ is isomorphic to $\mathbb{P}^1$, and the pullback of $\mathcal{O}_{\mathbb{P}^1}(1)$ to $E$ is an ample line bundle of degree $2|Q_0|$. If $Q_0$ is non-trivial, then this line bundle is very ample, and we therefore obtain a Galois embedding of $E$. By the previous discussion, we have that $E^d$ possesses a Galois embedding with Galois group $(Q_0\rtimes \mathbb{Z}/2\mathbb{Z})^d\rtimes S_d$, which is of order $2^dd!|Q_0|^d$. \qed

\subsection{Second proof of part (\ref{2})} Let $Q_0$ be a finite subgroup of $E$, let $G_0=Q_0^d\subseteq E^d$, and set $A=E^d$. We can see $A$ as 
$$A=\{(x_1,\ldots,x_{d+1})\in E^{d+1}:x_1+\cdots+x_{d+1}=0\};$$
and therefore we obtain an action of $S_{d+1}$ on $A$ by permutation. Notice that on $A$ itself the action is generated by the transformations
$$\sigma_1:(x_1,\ldots,x_d)\mapsto(x_2,x_1,x_3,\ldots,x_d)$$
$$\sigma_2:(x_1,\ldots,x_d)\mapsto(-(x_1+\ldots+x_d),x_1,\ldots,x_{d-1}).$$

Define $G=G_0\rtimes S_{d+1}$. As above we first take the quotient of $E^d$ by $G_0$ and obtain the abelian variety
$$B=(E/Q_0)^d.$$
The group $S_{d+1}$ acts on $B$ via the same permutations. Moreover, since we can think of $B$ as 
$$B=\{(y_1,\ldots,y_{d+1})\in (E/Q_0)^{d+1}:y_1+\cdots+y_{d+1}=0\},$$
we have that
\begin{eqnarray}\nonumber B/S_{d+1}&\simeq&\{p_1+\cdots+p_{d+1}\in\mbox{Sym}^{d+1}(E/Q_0):p_1+\cdots+p_{d+1}=0\mbox{ in }E/Q_0\}\\
\nonumber &\simeq&\{D\in\mbox{Div}(E/Q_0):D\geq0, D\sim_{\small\mbox{lin}} (d+1)[0]\}\\
\nonumber&=&|(d+1)[0]|\\
\nonumber&\simeq&\mathbb{P}^d
\end{eqnarray} 
where $[0]$ is the degree 1 divisor given by the origin of $E/Q_0$ and $\sim_{\small\mbox{lin}} $ stands for linear equivalence. Therefore, we have that $A/G\simeq\mathbb{P}^d$.

Let $\pi:A\to\mathbb{P}^d$ and $\psi:A\to B$ be the natural projections. We see that $\mathcal{O}_{\mathbb{P}^d}(1)$ is given by the hyperplane 
$$H=\{D'+[0]\in\mbox{Div}(E/Q_0):D'\geq0,D'\sim_{\small\mbox{lin}}  d[0]\},$$
and so $\mathcal{L}:=\pi^*\mathcal{O}_{\mathbb{P}^d}(1)\simeq\psi^*\mathcal{O}_B(\Theta_B+\ker(a))$, where $\Theta_B$ is the natural product polarization on $B$ and $a:B\to E/Q_0$ is the addition map.

If $Q_0$ is non-trivial and $d\geq2$ (or $|Q_0|\geq3$ and $d=1$), then $\mathcal{L}$ is very ample since it is the product of $|Q_0|^d$ ample line bundles, and by \cite[Theorem 1.1]{Bauer} the product of three ample line bundles is very ample.

Notice that the Galois group in this case is $G_0\rtimes S_{d+1}$ which is of order $(d+1)!|Q_0|^d$. \qed

\end{document}